\documentclass[11pt,reqno]{amsart}

\newcommand\bS{\mathbb{S}}

\newcommand{\nliminf}%
{\operatornamewithlimits{\underline{lim}\,}}
\newcommand{\nlimsup}%
{\operatornamewithlimits{\overline{lim}\,}}

\title[]{A simple proof of a result of A. Novikov}
\author[]{N.V. Krylov} 
\address{127 Vincent Hall, University of Minnesota,
 Minneapolis, MN, 55455}
\thanks{The work was partially supported
by an NSF grant}
\email{krylov@math.umn.edu}
\date{}%

\begin{document}

\begin{abstract}
 
We give simple proofs that for a continuous 
local martingale $M_{t}$ 
$$
\nliminf_{\varepsilon\downarrow0}\varepsilon
\log Ee^{(1-\varepsilon)
\langle M\rangle_{\infty}/2}<\infty
\Longrightarrow E\exp(M_{\infty}-\langle
 M\rangle_{\infty}/2)=1,
$$
$$
\nliminf_{\varepsilon\downarrow0}\varepsilon
\log\sup_{t\geq0}
Ee^{(1-\varepsilon)
M_{t}/2}<\infty\Longrightarrow 
E\exp(M_{\infty}-\langle M\rangle_
{\infty}/2)=1.
$$
\end{abstract}

\maketitle

1. Let $(\Omega,\mathcal{F},P)$ be a complete 
probability space and let
$M_{t}$ be a continuous local martingale 
on $(\Omega,\mathcal{F},P)$
such that $\langle M\rangle=\langle M
\rangle_{\infty}<\infty$ (a.s.).
Define
$$
M=M_{\infty},\,\,\,\rho=\rho(M)=e^{M-
\langle M\rangle/2},\,\,\,\rho_{t}=\rho_{t}
(M)=e^{M_{t}-
\langle M\rangle_{t}/2}.
$$
We will be discussing generalizations of the following 
 result of A.~Novikov (1973):
$$
Ee^{\langle M\rangle/2}<\infty\Longrightarrow 
E\rho=1.
$$
This result is quite important in many applications
related to absolute continuous change of probability measure
and, in particular, makes available Girsanov's theorem.

The original proof and other known proofs are 
rather complicated,
and here we want first to present an elementary
 proof of a somewhat
stronger result
\begin{equation}
                                           \label{12.9.1}
\lim_{\varepsilon\downarrow0}\varepsilon\log 
Ee^{(1-\varepsilon)
\langle M\rangle/2}=0\Longrightarrow E\rho=1.
\end{equation}
Then in n.~4  we show that $=0$ in (\ref{12.9.1})
 can be replaced with
$<\infty$.

It turns out that to prove (\ref{12.9.1}) 
it only suffices  to use the following two facts: 
\begin{equation}
                              \label{12.9.2}
 E\rho\leq1;\,\,\,                                                                 
 \exists\varepsilon>0:Ee^{(1+\varepsilon)
\langle M\rangle/2}<\infty\Longrightarrow E\rho=1.
\end{equation}
Indeed, if we accept (\ref{12.9.2}),
then under the condition in (\ref{12.9.1})
$$
Ee^{(1+\varepsilon)^{2}\langle (1-\varepsilon)M\rangle/2}
= Ee^{(1-\varepsilon^{2})^{2}\langle M\rangle/2}<\infty,
$$
which by (\ref{12.9.2}) and by   H\"older's 
inequality implies that
$$
1=E\rho((1-\varepsilon)M)=Ee^{(1-\varepsilon)
(M-\langle M\rangle/2)}
e^{(1-\varepsilon)\varepsilon\langle M\rangle/2}
$$
$$
\leq(E\rho)^{1-\varepsilon}
(Ee^{(1-\varepsilon)\langle M\rangle/2})^
{\varepsilon}.
$$
By letting $\varepsilon\downarrow0$ we get $1
\leq E\rho$, which
together with the first relation in (\ref{12.9.2})
 implies
our statement (\ref{12.9.1}).

2. In the same way we can improve a result of 
N. Kazamaki (1977). Let $\bS$ be the set of submartingales.
We claim that
\begin{equation}
                                  \label{12.9.3}
M^{+}_{\cdot}\in\bS,\quad
\lim_{\varepsilon\downarrow0}\varepsilon
\log\sup_{t\geq0}
Ee^{(1-\varepsilon)
M_{t}/2}=0\Longrightarrow E\rho=1.
\end{equation}

Here we use that
\begin{equation}
                               \label{12.9.4}
M^{+}_{\cdot}\in\bS,\quad
\exists\varepsilon>0:\,N_{0}:=\sup_{t\geq0}
Ee^{(1+\varepsilon)
M_{t}/2}<\infty\Longrightarrow E\rho=1.
\end{equation}
Then under the second condition in (\ref{12.9.3}) 
we have
$$
\sup_{t\geq0}Ee^{(1+\varepsilon)\{(1-\varepsilon)M_{t}\}/2}
=\sup_{t\geq0}Ee^{(1-\varepsilon^{2})M_{t}/2}<\infty,
$$
$$
1=E\rho((1-\varepsilon)M)=Ee^{(1-\varepsilon)
^{2}(M-\langle M\rangle/2)}
e^{(1-\varepsilon)\varepsilon M}
$$
$$
\leq(E\rho)^{(1-\varepsilon)^{2}}
(Ee^{(1-\varepsilon/(2-\varepsilon))  M /2})^
{\varepsilon(2-\varepsilon)},
$$
and we conclude as before. Note that the first condition in
\eqref{12.9.3} or 
\eqref{12.9.4} is satisfied if, for instance,
$E \langle M\rangle^{1/2}<\infty$ since then by Davis's inequality
$E\sup_{t\geq0}M_{t}^{+}<\infty$ and the local submartingale
$M^{+}_{t}$ is a submartingale.

3. Assertions (\ref{12.9.1}) and (\ref{12.9.3})
 are stronger than the
corresponding results of A.~Novikov and N.~Kazamaki. 
To show this for
(\ref{12.9.1}), take a one-dimensional Wiener 
process $w_{t}$ and let $\tau$
be the first exit time of $w_{t}$
from $(-\pi,\pi)$. Define $M_{t}=w_{t\wedge\tau}/4$.
 Then one can easily see that for 
$\varepsilon\downarrow0$
$$
[Ee^{(1-\varepsilon)\langle M\rangle/2)}
]^{\varepsilon}=
[Ee^{(1-\varepsilon)\tau/8}]^{\varepsilon}
=\big[\cos \frac{\sqrt{1-\varepsilon}}{2}\pi
\big]^{-\varepsilon}\rightarrow1,
$$
so that the assumption in (\ref{12.9.1}) is 
satisfied, whereas
$E\exp(\langle M\rangle/2)=\infty$ and  
 Novikov's criterion is not
applicable.

In the case of (\ref{12.9.3}) take $\tau$ 
to be an exponentially
distributed random variable independent of $w$.
 Specifically,
let $P(\tau>t)=e^{-t/2}$ and define 
$M=2w_{t\wedge\tau}$. Then
\begin{multline*}
[\sup_{t\geq0}Ee^{(1-\varepsilon)M_{t}/2}
]^{\varepsilon}
\\
=
\left[\sup_{t\geq0}\frac{1}{2}\left\{\int_{0}^{t}e^{-s/2}
Ee^{(1-\varepsilon)w_{s}}\,ds+
Ee^{(1-\varepsilon)w_{t}}\int_{t}^{\infty}
e^{-s/2}\,ds\right\}\right]^{\varepsilon}
\\
=\left[\frac{1}{2}\int_{0}^{\infty}
e^{-s(1-(1-\varepsilon)^{2})/2}\,ds
\right]^{\varepsilon}\rightarrow1.
\end{multline*}
At the same time $\sup_{t\geq0}E
\exp(M_{t}/2)=\infty$,
and  Kazamaki's criterion is not applicable.

4. Now we show further 
improvements of (\ref{12.9.1})
and (\ref{12.9.3}):
\begin{equation}
                                           \label{3.19.1}
\nliminf_{\varepsilon\downarrow0}\varepsilon\log 
Ee^{(1-\varepsilon)
\langle M\rangle/2}<\infty\Longrightarrow E\rho=1,
\end{equation}
\begin{equation}
                                            \label{3.19.2}
M^{+}_{\cdot}\in\bS,
\quad\nliminf_{\varepsilon\downarrow0}\varepsilon
\log\sup_{t\geq0}
Ee^{(1-\varepsilon)
M_{t}/2}<\infty\Longrightarrow E\rho=1.
\end{equation}

To prove (\ref{3.19.1}) we proceed as in the proof of (\ref{12.9.1})
and  we write 
$$
1=E\rho((1-\varepsilon)M)=Ee^{(1-\varepsilon)
(M-\langle M\rangle/2)}
e^{(1-\varepsilon)\varepsilon\langle
 M\rangle/2}I_{\langle M\rangle\leq T}
$$
$$
+EI_{\langle M\rangle> T}e^{(1-\varepsilon)
(M-\langle M\rangle/2)}
e^{(1-\varepsilon)\varepsilon\langle M\rangle/2}
$$
$$
\leq(E\rho)^{1-\varepsilon}
(Ee^{(1-\varepsilon)\langle M\rangle/2}I_{\langle M\rangle\leq T})^
{\varepsilon}+
(E\rho I_{\langle M\rangle> T})^{1-\varepsilon}
(Ee^{(1-\varepsilon)\langle M\rangle/2})^
{\varepsilon},
$$
where $T$ is a constant, $T\in(0,\infty)$.
As $\varepsilon\downarrow0$, we get $$1
\leq E\rho+{\rm const} \,
E\rho I_{\langle M\rangle> T},$$  which gives $1\leq E\rho$
after letting $T\rightarrow\infty$. In like manner
(\ref{3.19.2}) is proved.

5. For the sake of completeness we also 
present the proofs of (\ref{12.9.2})
and (\ref{12.9.4}). The first relation 
in (\ref{12.9.2}) is true because
$\rho_{t}$ is a nonnegative local martingale 
(by It\^o's formula). From this
we get
$$
Ee^{M_{t}/2}=Ee^{M_{t}/2-\langle M\rangle_{t}/4}
e^{\langle M\rangle_{t}/4}
\leq (Ee^{\langle M\rangle_{t}/2})^{1/2}\leq 
(Ee^{\langle M\rangle/2})^{1/2},
$$
and it remains only to prove (\ref{12.9.4}).

To do the last step
take $\kappa>1$ and $p>1$, define 
$\gamma=(p\kappa)^{-1/2}$
and $q=p/(p-1)$, and let $\tau_{n}$ be a sequence
of bounded stopping times localizing $\rho$.
Observe that by Doob's inequality for moments 
of the martingales
$\rho_{t\wedge\tau_{n}}$ and the assumption that
$M^{+}_{t}$ is a submartingale we have
$$
E\sup_{t\leq\tau_{n}}\rho_{t}^{\kappa}\leq
N_{1} E\rho_{\tau_{n}}^{\kappa}=
N_{1} Ee^{\gamma\kappa M_{\tau_{n}}-
\kappa\langle M_{\tau_{n}}\rangle/2}e^{(1-\gamma)
\kappa M_{\tau_{n}}}
$$
\begin{equation}
                                                  \label{4.29.1}
\leq N_{1}(E\rho_{\tau_{n}}(p\gamma\kappa M))^{1/p}
\sup_{t>0}(Ee^{\delta M_{t}^{+}})^{1/q}\leq N_{1}
(N_{0}+1)^{1/q},
\end{equation}
where $\delta=(1-\gamma)\kappa q$, $N_{1}$ 
is a constant and the last inequality is true
if $\delta\leq(1+\varepsilon)/2$. The latter is easy to accommodate
since $\delta\to1/2$ as first $\kappa\downarrow1$ and then 
$ p\downarrow1$.  Therefore, given that 
the condition in
(\ref{12.9.4}) is satisfied, we can find 
$\kappa>1$ and $p>1$ such that
$\delta\leq(1+\varepsilon)/2$, and then
by  \eqref{4.29.1} and Fatou's lemma
$E\sup_{t\geq0}\rho_{t}^{\kappa}<\infty$ and 
$E\sup_{t\geq 0}\rho_{t}<\infty$. Finally, by using the
dominated convergence theorem we conclude
$$
E\rho=E\lim_{n\rightarrow\infty}\rho_{\tau_{n}}=
\lim_{n\rightarrow\infty}E\rho_{\tau_{n}}=1.
$$

\end{document}